\documentclass[12pt,a4paper,reqno]{amsart}
\usepackage[latin1]{inputenc}
\usepackage{amsmath}
\usepackage{amsfonts}
\usepackage{amssymb}
\usepackage{amsmath,bbm,amssymb,amsxtra}
\usepackage{enumerate}
\allowdisplaybreaks
\usepackage{epsfig}
\usepackage{ae}
    \def\Xint#1{\mathchoice
    {\XXint\displaystyle\textstyle{#1}}%
    {\XXint\textstyle\scriptstyle{#1}}%
    {\XXint\scriptstyle\scriptscriptstyle{#1}}%
    {\XXint\scriptscriptstyle\scriptscriptstyle{#1}}%
    \!\int}
    \def\XXint#1#2#3{{\setbox0=\hbox{$#1{#2#3}{\int}$ }
    \vcenter{\hbox{$#2#3$ }}\kern-.6\wd0}}
    \def\dashint{\Xint-}

\theoremstyle{definition}

\newtheorem{lemma}{Lemma}[section]

\newtheorem{proposition}[lemma]{Proposition}

\newtheorem{theorem}[lemma]{Theorem}

\newtheorem*{acknowledgements}{Acknowledgements}
\newtheorem{definition}[lemma]{Definition}

\newcommand{\prop}[1]{\begin{proposition}\label{#1}
\sl }
\newcommand{\eprop}{\end{proposition}}
\newcommand{\thm}[1]{\begin{theorem}\label{#1}
\sl }
\newcommand{\ethm}{\end{theorem}}

\newcommand{\lem}[1]{\begin{lemma}\label{#1}
\sl }
\newcommand{\elem}{\end{lemma}}

\newcommand{\defin}[1]{\begin{definition}\label{#1}
\sl }
\newcommand{\edefin}{\end{definition}}

\newcommand{\beqno}{\begin{eqnarray*}}
\newcommand{\eeqno}{\end{eqnarray*}}
\newcommand{\beqla}[1] {\begin {eqnarray}\label{#1}}
\def\eeq {\end {eqnarray}}
\newcommand{\beq}{\begin {eqnarray}}

\newcommand{\real}{{\mathbb R}}

\newcommand{\integer}{{\mathbb Z}}
\newcommand{\nanu}{{\mathbb N}}
\newcommand{\complex}{{\mathbb C}}

\newcommand{\sh}{{\mathcal S}}
\newcommand{\sdistr}{{\mathcal S}'}
\newcommand{\poly}{{\mathcal P}}

\newcommand{\supp}{{\rm supp}\,}

\newcommand{\D}{\mathbb{D}}
\newcommand{\hlmax}{\mathcal{M}}
\newcommand{\mca}{\mathcal{A}}
\newcommand{\df}{\dot{F}}
\newcommand{\dfp}{\dot{F}^{s,\tau}_{p,q}(\real^n)}
\newcommand{\dfpp}{\dot{F}^s_{p,q}(\real^n)}
\newcommand{\db}{\dot{B}}
\newcommand{\dbp}{\dot{B}^{s,\tau}_{p,q}(\real^n)}
\newcommand{\dbpp}{\dot{B}^s_{p,q}(\real^n)}

\newcommand{\dmp}{\dot{M}^{s,\tau}_{p,q}(\real^n)}
\newcommand{\dmpp}{\dot{M}^s_{p,q}(\real^n)}

\newcommand{\dnp}{\dot{N}^{s,\tau}_{p,q}(\real^n)}

\newcommand{\dnpp}{\dot{N}^s_{p,q}(\real^n)}

\title[Characterizations of Besov-type and Triebel-Lizorkin-type spaces]{Pointwise and grand maximal function characterizations of Besov-type and Triebel-Lizorkin-type spaces}

\author[Soto]{Tom\'as Soto}
\address{Department of Mathematics and Statistics, University of
Helsinki, PO~Box~68, FI-00014 Helsinki, Finland}
\email{tomas.soto@helsinki.fi}

\thanks{The author was supported by the Finnish CoE in Analysis and Dynamics Research}

\keywords{Besov-type space, function space, grand maximal function, Haj\l asz gradient, Triebel-Lizorkin-type space}

\subjclass[2000]{Primary: 42B35; Secondary: 46E35}

\begin{document}

\maketitle

\begin{abstract}
In this note, we establish characterizations for the homogeneous Besov-type spaces $\dot{B}^{s,\tau}_{p,q}(\mathbb{R}^n)$ and Triebel-Lizorkin-type spaces $\dot{F}^{s,\tau}_{p,q}(\mathbb{R}^n)$, introduced by Yang and Yuan, through fractional Haj\l asz-type gradients for suitable values of the parameters $p$, $q$ and $\tau$ when $0 < s < 1$, and through grand Littlewood-Paley-type maximal functions for all admissible values of the parameters. These characterizations extend the characterizations obtained by Koskela, Yang and Zhou for the standard homogeneous Besov and Triebel-Lizorkin spaces.
\end{abstract}

\section{Introduction}\label{se:introduction}

The main purpose of this note is to establish pointwise characterizations of the Besov-type and Triebel-Lizorkin-type function spaces, introduced by Yang and Yuan in \cite{YY} and \cite{YY2}, through Haj\l asz-type gradients. Characterizations of this type go back to Haj\l asz's pointwise characterization of the classical Sobolev spaces \cite{H}, and they have found many applications in both the Euclidean setting as well as in the setting of more general metric measure spaces. The families of function spaces considered in this paper include the standard Besov and Triebel-Lizorkin spaces as well as the fractional Morrey-Sobolev spaces for smoothness indices $s \in (0,1)$ as special cases.

To begin with, we first recall the definitions of the standard (homogeneous) Triebel-Lizorkin and Besov spaces. For a dimension $n \in \nanu := \{1,2,3,\cdots\}$, which shall be fixed throughout the paper, we denote by $\sh(\real^n)$ the class of Schwartz functions, i.e.~the class of complex-valued $C^{\infty}(\real^n)$ functions $\phi$ for which
\[
  \|\phi\|_{\sh_{k,m}} := \sup_{|\gamma| \leq k,\; x\in\real^n} (1 + |x|)^m |\partial^\gamma \phi(x)|
\]
is finite for all $k$, $m \in \nanu_0 := \nanu \cup \{0\}$; here $|\gamma| = |\gamma_1|+\cdots+|\gamma_n|$ and $\partial^\gamma = \partial^{\gamma_1}_{x_1}\cdots\partial^{\gamma_n}_{x_n}$ for all multi-indices $\gamma = (\gamma_1,\cdots,\gamma_n) \in \nanu_0^n$. The seminorms $\|\cdot\|_{\sh_{k,m}}$ induce a locally convex topology on $\sh(\real^n)$. We denote by $\sdistr(\real^n)$ the class of tempered distributions, i.e.~the class of continuous complex-valued linear functionals on $\sh(\real^n)$, and equip $\sdistr(\real^n)$ with the topology induced by the mappings $f \mapsto \langle f,\phi \rangle$, $\phi \in \sh(\real^n)$. For standard facts about the Schwartz space and tempered distributions, particularly their Fourier-transforms and convolutions, we refer to e.g.~\cite{G}.

Following \cite{FJ}, $\varphi$ and $\psi$ shall throughout the paper be fixed elements of $\sh(\real^n)$ satisfying
\[
  \supp\widehat{\varphi},\; \supp\widehat{\psi} \subset \big\{\xi \in \real^n : \ 2^{-1} \leq |\xi| \leq 2\big\},
\]
\[
  |\widehat{\varphi}(\xi)|,\, |\widehat{\psi}(\xi)| \geq c > 0 \;\, \text{when } 3/5 \leq |\xi| \leq 5/3
\]
and
\[
  \sum_{j\in\integer} \overline{\widehat{\varphi}(2^j\xi)}\widehat{\psi}(2^j\xi) = 1 \;\, \text{when } \xi \neq 0.
\]
For $\phi \in \sh(\real^n)$ and $j \in \integer$, we write $\phi_j$ for the Schwartz function $x \mapsto 2^{jn}\phi(2^j x)$.

For $s \in \real$, $0 < p < \infty$ and $0 < q \leq \infty$, the homogeneous Triebel-Lizorkin space $\dfpp$ is defined as the class of tempered distributions $f$ for which
\[
  \|f\|_{\dfpp} = \left(\int_{\real^n} \bigg[\sum_{j\in\integer} \Big(2^{js}|\varphi_j*f(z)|\Big)^q\bigg]^{p/q} dz\right)^{1/p}
\]
is finite, with the obvious modification made when $q = \infty$. For $s \in \real$, $0 < p \leq \infty$ and $0 < q \leq \infty$, the homogeneous Besov space $\dbpp$ is defined as the class of tempered distributions $f$ for which
\[
  \|f\|_{\dbpp} = \left(\sum_{j\in\integer}\Big(2^{js} \|\varphi_j*f\|_{L^p(\real^n)}\Big)^q\right)^{1/q}
\]
is finite, with again the obvious modification made when $q = \infty$. It is well known that after quotienting out the tempered distributions whose Fourier-transforms are supported at the origin, i.e.~the polynomials, $\dfpp$ and $\dbpp$ become quasi-Banach spaces, independent of the choice of $\varphi$ in the sense that two admissible choices induce equivalent quasinorms; see for instance \cite{T}.

The following Triebel-Lizorkin-type and Besov-type spaces were introduced by Yang and Yuan in \cite{YY} and \cite{YY2}. For $s \in \real$, $0 < p < \infty$, $0 < q \leq \infty$ and $0 \leq \tau < \infty$, the homogeneous Triebel-Lizorkin-type space $\dfp$ is defined as the class of tempered distributions $f$ for which
\[
  \|f\|_{\dfp} = \sup_{x\in\real^n,\;\ell \in \integer} \frac{1}{|B(x,2^{-\ell})|^{\tau} }\left(\int_{B(x,2^{-\ell})} \bigg[ \sum_{j \geq \ell}\Big(2^{js} |\varphi_j*f(z)|\Big)^q \bigg]^{p/q} dz\right)^{1/p}
\]
is finite, with the obvious modification when $q = \infty$. For $s \in \real$, $0 < p \leq \infty$, $0 < q \leq \infty$ and $0 \leq \tau < \infty$, the homogeneous Besov-type space $\dbp$ is defined as the class of tempered distributions $f$ for which
\[
  \|f\|_{\dbp} = \sup_{x\in\real^n,\;\ell \in \integer} \frac{1}{|B(x,2^{-\ell})|^{\tau}} \left( \sum_{j \geq \ell}\Big(2^{js} \|\varphi_j * f\|_{L^p(B(x,2^{-\ell}))}\Big)^q\right)^{1/q}
\]
is finite, with again the obvious modification when $q = \infty$. Again, these spaces become quasinormed spaces after quotienting out the polynomials. Actually, in the definitions of \cite{YY} and \cite{YY2}, the supremum is taken over dyadic cubes instead of balls with dyadic radii, but it is quite easy to see that the above definitions yield equivalent quasinorms. It is also known that these spaces are independent of the choice of $\varphi$ (\cite[Corollary 3.1]{YY2}) and that they are quasi-Banach spaces (\cite[Proposition 2.2]{SYY} and the references therein). Here are a couple examples of how these spaces coincide\footnote{Here and in the sequel, $X = Y$ for function spaces $X$ and $Y$ means that they embed continuously into each other.} with other well-known function spaces:
\begin{itemize}
\item $\df^{s,0}_{p,q}(\real^n) = \dfpp$ and $\db^{s,0}_{p,q}(\real^n) = \dbpp$ for all admissible parameters.
\item $\df^{s,1/p}_{p,q}(\real^n) = \df^{s}_{\infty,q}(\real^n)$ for all admissible parameters; in particular $\df^{0,1/p}_{p,2}(\real^n) = BMO(\real^n)$. 
\item $\db^{s,\frac1u - \frac1p}_{u,\infty}(\real^n) = \dot{\mathcal{N}}^{s}_{p,u,\infty}(\real^n)$ for $0 < u \leq p < \infty$ and $s\in\real$, where $\dot{\mathcal{N}}^{s}_{p,u,q}(\real^n)$ stands for the homogeneous Besov-Morrey space, i.e.~the Besov-type space based on the Morrey space $M^p_u(\real^n)$ instead of $L^p(\real^n)$, introduced in \cite{KY} and \cite{M}.
\item $\df^{s,\frac1u - \frac1p}_{u,q}(\real^n) = \dot{\mathcal{E}}^{s}_{p,u,q}(\real^n)$ for $0 < u \leq p < \infty$, $0 < q \leq \infty$ and $s \in \real$, where $\dot{\mathcal{E}}^{s}_{p,u,q}(\real^n)$ stands for the homogeneous Triebel-Lizorkin-Morrey space, i.e.~the Triebel-Lizorkin-type space based on the Morrey space $M^p_u(\real^n)$ instead of $L^p(\real^n)$.
\item For $0 < \alpha < \min(1,\frac{n}{2})$, the space $\df^{\alpha,\frac12 - \frac{\alpha}{n}}_{2,2}(\real^n)$ coincides with the space $Q_\alpha(\real^n)$ introduced in \cite{EJPX}. 
\item For $\frac1p < \tau < \infty$ and all admissible values of $p$, $q$ and $s$,
\beqla{eq:zygmund}
  \dfp = \dot{F}^{s + n(\tau - \frac1p)}_{\infty,\infty}(\real^n) \quad {\rm and} \quad \dbp = \dot{B}^{s + n(\tau - \frac1p)}_{\infty,\infty}(\real^n).
\eeq
Recall that for $s > 0$, $\dot{F}^{s}_{\infty,\infty}(\real^n)$ and $\dot{B}^{s}_{\infty,\infty}(\real^n)$ both coincide with the homogeneous H\"older-Zygmund space of order $s$; see e.g.~\cite[Theorem 6.3.6]{G} and \cite[Section 5]{FJ}.
\end{itemize}
We refer to e.g.~\cite{SYY} and \cite{YY4} for the definitions of the spaces $\dot{F}^s_{\infty,q}(\real^n)$, $BMO(\real^n)$, $\dot{\mathcal{E}}^{s}_{p,u,q}(\real^n)$, $\dot{\mathcal{N}}^{s}_{p,u,q}(\real^n)$ and $Q_\alpha(\real^n)$ (which we shall not need in the sequel) and to \cite[Proposition 1 and Theorem 2]{YY4} as well as the references therein for details about the above coincidences. We also refer to \cite{SYY} for a detailed discussion of the history of the spaces in question.

Inspired by \cite{KYZ2}, we now define function spaces analogous to $\dfp$ and $\dbp$ through Haj\l asz-type gradients. If $u : \real^n \to \complex$ is a (Lebesgue-)measurable function and $0 < s < \infty$, $\D^s(u)$ stands for the class of all \emph{fractional s-Haj\l asz gradients} of $u$, i.e.~the class of sequences $\overrightarrow{g} = (g_k)_{k\in\integer}$ of measurable functions $g_k : \real^n \to [0,\infty]$ for which there exists a set $E \subset \real^n$ of measure zero such that, for all $k \in \integer$,
\[
  |u(x) - u(y)| \leq |x-y|^s\left(g_k(x) + g_k(y)\right)
\]
whenever $x$, $y \in \real^n \backslash E$ and $2^{-k-1} \leq |x-y| < 2^{-k}$.

For $0 < p < \infty$, $0 < q \leq \infty$, $0 < s < \infty$ and $0 \leq \tau < \infty$, we define $\dmp$ as the class of measurable functions $u$ such that
\[
  \|u\|_{\dmp} := \inf_{\overrightarrow{g} \in\D^s(u)} \sup_{x\in\real^n,\;\ell\in\integer} \frac{1}{|B(x,2^{-\ell})|^{\tau}}\left(\int_{B(x,2^{-\ell})} \bigg( \sum_{k \geq \ell}g_k(z)^q \bigg)^{p/q} dz\right)^{1/p}
\]
is finite, and for $0 < p \leq \infty$, $0 < q \leq \infty$, $0 < s < \infty$ and $0 \leq \tau < \infty$, we define $\dnp$ as the class of measurable functions $u$ such that
\[
  \|u\|_{\dnp} := \inf_{\overrightarrow{g} \in\D^s(u)} \sup_{x\in\real^n,\;\ell\in\integer} \frac{1}{|B(x,2^{-\ell})|^{\tau}} \left( \sum_{k \geq \ell}\|g_k\|_{L^p(B(x,2^{-\ell}))}^q\right)^{1/q}
\]
is finite, with the obvious modifications for $q = \infty$ in both cases. $\dmp$ and $\dnp$ become quasinormed spaces after quotienting out the functions that are constant almost everywhere. For $\tau = 0$, they coincide with the spaces $\dmpp$ and $\dnpp$ introduced in \cite{KYZ2}. Note that analogous function spaces could well be defined on any metric measure space instead of just $\real^n$; see \cite{KYZ2} and \cite{GKZ} for the case $\tau = 0$. 

Our main result is the following pointwise characterization of the elements of $\dfp$ and $\dbp$ for $s \in (0,1)$ which generalizes the result obtained for $\tau = 0$ in \cite[Theorem 3.2]{KYZ2}.

\thm{th:main}
(i) For $s \in (0,1)$, $p \in (\frac{n}{n+s},\infty)$, $q \in (\frac{n}{n+s},\infty]$ and $\tau \in [0,\frac{1}{p} + \frac{1-s}{n})$, we have $\dmp = \dfp$ with equivalent quasinorms.

(ii) For $s \in (0,1)$, $p \in (\frac{n}{n+s},\infty]$, $q \in (0,\infty]$ and $\tau \in [0,\frac{1}{p} + \frac{1-s}{n})$, we have $\dnp = \dbp$ with equivalent quasinorms.
\ethm

{\noindent Note that if $0 < s < 1$ and $\tau \geq 1/p + (1-s)/n$, then $s + n(\tau - 1/p) \geq 1$, so in view of \eqref{eq:zygmund} we do not expect that the range of $\tau$ in the theorem could be improved. The theorem follows from Theorems \ref{th:grand} and \ref{th:pointwise} below, which are of independent interest and whose setting we shall explain next.}

Inspired by \cite{KYZ} and \cite{KYZ2}, we define ``grand'' counterparts of $\dfp$ and $\dbp$. For all $N \in \nanu_0 \cup \{-1\}$, and $m$, $l \in \nanu_0$, let
\[
  \mca^{l}_{N,m} = \left\{\phi \in \sh(\real^n) : \int_{\real^n} x^\gamma \phi(x) dx = 0 \text{ when } |\gamma| \leq N \text{ and } \|\phi\|_{\sh_{N+l+1,m}} \leq 1 \right\},
\]
where the moment condition is interpreted to be void when $N = -1$. For $s \in \real$, $0 < p < \infty$, $0 < q \leq \infty$, $0 \leq \tau < \infty$ and $N$, $m$ and $l$ as above, we define the grand Triebel-Lizorkin-type space $\mca^l_{N,m}\dfp$ as the class of tempered distributions $f$ for which
\[
  \|f\|_{\mca^l_{N,m}\dfp} := \sup_{x\in\real^n,\;\ell\in\integer} \frac{1}{|B(x,2^{-\ell})|^{\tau}}\left(\int_{B(x,2^{-\ell})} \bigg[ \sum_{j \geq \ell}\Big(2^{js} \sup_{\phi \in \mca^l_{N,m}}|\phi_j*f(z)|\Big)^q \bigg]^{\frac{p}{q}}dz\right)^{\frac1p}
\]
is finite, with the obvious modification for $q = \infty$. For $s \in \real$, $0 < p \leq \infty$, $0 < q \leq \infty$, $0 \leq \tau < \infty$ and $N$, $m$ and $l$ as above, we define the grand Besov-type space $\mca^l_{N,m}\dbp$ as the class of tempered distributions $f$ for which
\[
  \|f\|_{\mca^l_{N,m}\dbp} := \sup_{x\in\real^n,\;\ell\in\integer} \frac{1}{|B(x,2^{-\ell})|^{\tau}} \left( \sum_{j \geq \ell}\bigg(2^{js} \Big\|\sup_{\phi \in \mca^l_{N,m}} |\phi_j * f| \Big\|_{L^p(B(x,2^{-\ell}))}\bigg)^q\right)^{\frac1q}
\]
is finite, where the supremum is taken pointwise and the obvious modification is made for $q = \infty$. The two quantities defined above are quasinorms when $N = -1$, and when $N \in \nanu_0$, they become quasinorms after quotienting out the polynomials with degree at most $N$. We obviously have $\mca^l_{N,m}\df^{s,0}_{p,q}(\real^n) = \mca^l_{N,m}\dfpp$ and $\mca^l_{N,m}\db^{s,0}_{p,q}(\real^n) = \mca^l_{N,m}\dbpp$ for all admissible parameters, where $\mca^l_{N,m}\dfpp$ and $\mca^l_{N,m}\dbpp$ are the grand Triebel-Lizorkin and grand Besov spaces defined in \cite{KYZ} and \cite{KYZ2}. The following result extends \cite[Theorem 3.1]{KYZ2}.

\thm{th:grand}
(i) Let $0 < p < \infty$, $0 < q \leq \infty$, $s \in \real$ and $J = n/\min(1,p,q)$. If the integers $N \geq -1$, $m \geq 0$ and $l \geq 0$ satisfy
\beqla{eq:grand-parameters}
  N+1 > \max(s,J-n-s) \quad {\rm and} \quad m > \max(J,n+N+1),
\eeq
then $\dfp = \mca^l_{N,m}\dfp$ for all $\tau \in [0, \frac1p]$.

(ii) Let $0 < p \leq \infty$, $0 < q \leq \infty$, $s \in \real$ and $J = n/\min(1,p)$. If the integers $N \geq -1$, $m \geq 0$ and $l \geq 0$ are related by \eqref{eq:grand-parameters}, then $\dbp = \mca^l_{N,m}\dbp$ for all $\tau \in [0, \frac1p]$.

(iii) The results of (i) and (ii) actually hold for all $\tau \in [0, \frac1p + \frac{\epsilon}{2n})$, where
\[ \epsilon = \min\big(2(N+1-s),\,m-J,\,2(N+1+n+s-J)\big) > 0. \]
\ethm
{\noindent The proof of the theorem is presented in section \ref{se:grand}.}
We also have the following result which generalizes the analogous results for $\tau = 0$ obtained in \cite[Theorem 3.2]{KYZ2}.

\thm{th:pointwise}
Suppose that $s \in (0,1)$ and $m \geq n+1$ or $s = 1$ and $m \geq n+2$.

(i) For $p \in (\frac{n}{n+s},\infty)$, $q \in (\frac{n}{n+s},\infty]$, $\tau \in [0,\frac{1}{p} + \frac{1-s}{n})$ and $l \in \nanu_0$, we have $\dmp = \mca^l_{0,m} \dfp$ with equivalent quasinorms.

(ii) For $p \in (\frac{n}{n+s},\infty]$, $q \in (0,\infty]$, $\tau \in [0,\frac{1}{p} + \frac{1-s}{n})$ and $l \in \nanu_0$, we have $\dnp = \mca^l_{0,m} \dbp$ with equivalent quasinorms.
\ethm

{\noindent The proof of this is presented in section \ref{se:pointwise}.}

Previously, the Besov-type and Triebel-Lizorkin type spaces have been characterized through difference integrals (see e.g.~\cite{D} where the non-homogeneous versions of the spaces are treated) and through Peetre-type maximal functions as well as local means in \cite{YY3}. We refer to \cite{SYY} for a variety of other characterizations.

We end this section with some notation conventions. We write $|A|$ for the $n$-dimensional Lebesgue measure of a measurable set $A$, and $\dashint_A f$ or $f_A$ for $|A|^{-1} \int_A f(x) dx$ whenever the latter quantity is well-defined. For two non-negative functions $f$ and $g$ with the same domain, we write $f \lesssim g$ if $f \leq Cg$ for some positive and finite constant $C$, usually independent of some paramters; $f \approx g$ means that $f \lesssim g$ and $g \lesssim f$. For real numbers $x$ and $y$, we may write $x\land y$ for $\min(x,y)$. For a dyadic cube $Q := 2^{-j}(k + [0,1]^n)$, $j \in \integer$, $k \in \integer^n$, we let $\ell(Q) = 2^{-j}$ and $x_Q = 2^{-j}k$. For a ball $B$ of $\real^n$, we denote by $\lambda B$, $\lambda > 0$, the ball cocentric with $B$ with radius $\lambda$ times the radius of $B$. 

We write $\poly(\real^n)$ for the vector space of polynomials in $\real^n$ and $\poly_N(\real^n)$ for the vector space of polynomials with degree at most $N$ when $N \in \nanu_0$; when $N < 0$, we let $\poly_N(\real^n) := \{0\}$. We shall frequently abuse notation by writing $\dfp$ for both the class of tempered distributions $f$ such that $\|f\|_{\dfp}$ is finite as well as the function space $\{ f \in \sdistr(\real^n)/\poly(\real^n) : \|f\|_{\dfp} < \infty\}$, and similarly for other spaces defined in this section as well. When talking about $\dfp$ and $\dbp$ in the same sentence with some indicated parameter range, it is understood that the possibility $p = \infty$ is excluded in the case of $\dfp$, and similarly for the other pairs of spaces defined in this section.

\section{Proof of Theorem \ref{th:grand}}\label{se:grand}

The argument below modifies the proof of \cite[Theorem 3.1]{KYZ2}; we have left out some details that carry over unchanged.

\begin{proof}[Proof of Theorem \ref{th:grand}]
(i) First, if a tempered distribution $f$ belongs to $\mca^l_{N,m}\dfp$, then $\|f\|_{\dfp} \lesssim \|f\|_{\mca^l_{N,m}\dfp} < \infty$ since our fixed $\varphi$ is a constant multiple of an element of $\mca^l_{N,m}$. In this sense we have $\mca^l_{N,m}\dfp \subset \dfp$.

Conversely, if a tempered distribution $f$ belongs to $\dfp$, it is shown in the proof of \cite[Lemma 4.2]{YY2} that
\[
  \|\partial^\gamma (\widetilde{\psi_j}*\varphi_j*f)\|_{L^{\infty}(\real^n)} \lesssim 2^{j(|\gamma|+n/p-s-n\tau)} \|f\|_{\dfp},
\]
for all integers $j$ and multi-indices $\gamma$, so the discussion in \cite[pp.~153--155]{FJ} applies: if $|\gamma| > s + n\tau - n/p$, then $\sum_{j \leq 1}\|\partial^\gamma (\widetilde{\psi_j}*\varphi_j*f)\|_{L^{\infty}(\real^n)} < \infty$, so there exists polynomials $(P_K)_{K\in \nanu}$ in $\poly_L(\real^n)$ with $L = \lfloor s + n\tau - n/p\rfloor$ and a polynomial $P_f$ such that
\[
  f + P_f = \lim_{K \to \infty} \left(\sum_{j = -K}^{\infty} \widetilde{\psi_j}*\varphi_j*f + P_K\right)
\]
with convergence in $\sdistr(\real^n)$; here the polynomial $P_f$ is unique modulo $\poly_L(\real^n)$ in the sense that if $\psi^{(i)}$, $\varphi^{(i)}$, $(P_K^{i})_{K\in\nanu}$ and $P^{i}_f$ are, for $i \in \{1,2\}$, two choices of admissible functions as above, then $P^{(1)}_f - P^{(2)}_f \in \poly_L(\real^n)$. Furthermore, $L < s + \frac{\epsilon}2 \leq N+1$, and it is easily checked that for $f_1,\,f_2 \in \dfp$, $(f_1 + P_{f_1}) - (f_2 + P_{f_2})$ is an element of $\poly_N(\real^n)$ if and only if $f_1 - f_2$ is a polynomial. The rule $f \mapsto f + P_f =: \widetilde{f}$ thus yields a well-defined, linear and injective mapping from the function space $\dfp$ into $\sdistr(\real^n)/\poly_N(\real^n)$. The plan now is to show that this mapping actually takes $\dfp$ continuously into $\mca^l_{N,m}\dfp$.

Since $m > n+N+1$, by the proof of \cite[Theorem 1.2]{KYZ} we have
\beqla{eq:ad-estimate}
  |\phi_j*\widetilde{f}(z)| \lesssim \sum_{\ell(Q) = 2^{-j}} \bigg(\sum_{R} a_{QR}|t_R|\bigg) |Q|^{-1/2} \chi_Q(z)
\eeq
for all $\phi \in\mca^l_{N,m}$, $z \in \real^n$ and $j \in \integer$, with the implied constant independent of these parameters. Here the outer sum is taken over all dyadic cubes $Q$ of $\real^n$ with side length $2^{-j}$ and the inner sum over all dyadic cubes $R$ of $\real^n$; $t_R$ denotes $\langle f,\widetilde{\psi_R} \rangle$ with $\psi_R(x) = 2^{-in/2}\psi_i(x - x_R)$ for all dyadic cubes $R$ with side length $2^{-i}$; and
\beqla{eq:ad-estimate-2}
  a_{QR} = 2^{-|i-j|(n/2 + N + 1)}\left(1 + 2^{\min(i,j)}|x_Q - x_R|\right)^{-m}
\eeq
for all cubes $Q$ and $R$ with side lengths $2^{-j}$ and $2^{-i}$ respectively. The cubes $Q$ in \eqref{eq:ad-estimate} are pairwise disjoint, so for any dyadic cube $P$ with side length $2^{-\ell}$ we have
\begin{align*}
\, &  \frac{1}{|P|^{\tau}} \left(\int_P \bigg[\sum_{j \geq \ell}\Big(2^{js}\sup_{\phi \in \mca}|\phi_j * \widetilde{f}(z)|\Big)^q \bigg]^{p/q}dz  \right)^{1/p} \\
\lesssim\;  & \frac{1}{|P|^{\tau}}\left(\int_P \bigg[\sum_{j \geq \ell} 2^{jsq}\sum_{\ell(Q) = 2^{-j}} \bigg(\sum_{R} a_{QR}|t_R|\bigg)^q |Q|^{-q/2} \chi_Q(z) \bigg]^{p/q}dz  \right)^{1/p} \\
\approx\; & \frac{1}{|P|^{\tau}}\left(\int_P \bigg[\sum_{Q \subset P} |Q|^{-sq/n-q/2}  \bigg(\sum_{R} a_{QR}|t_R|\bigg)^q \chi_Q(z) \bigg]^{p/q}dz  \right)^{1/p} \\
\lesssim\; & \bigg\|\bigg(\sum_{R} a_{QR}|t_R|\bigg)_Q\bigg\|_{\dot{f}^{s,\tau}_{p,q}(\real^n)},
\end{align*}
where $\dot{f}^{s,\tau}_{p,q}(\real^n)$ stands for the space of sequences $(b_Q)_Q$ of complex numbers, indexed by the dyadic cubes of $\real^n$, such that
\[
  \|(b_Q)_Q\|_{\dot{f}^{s,\tau}_{p,q}(\real^n)} = \sup_{P} \frac{1}{|P|^{\tau}}\left(\int_P \bigg[\sum_{Q \subset P} \Big(|Q|^{-s/n-1/2} |b_Q| \chi_Q(x)\Big)^q \bigg]^{p/q}dx  \right)^{1/p},
\]
where the supremum is taken over all dyadic cubes $P$ of $\real^n$, is finite. From \eqref{eq:ad-estimate-2} it is easy to check that for $\epsilon > 0$ as in the statement of part (iii), we have
\[
  a_{QR} \leq \left(\frac{\ell(Q)}{\ell(R)}\right)^s \left(1 + \frac{|x_Q-x_R|}{\max(\ell(Q),\ell(R))}\right)^{-J-\epsilon}\min\left[ \left(\frac{\ell(Q)}{\ell(R)}\right)^{\frac{n+\epsilon}{2}},\left(\frac{\ell(R)}{\ell(Q)}\right)^{J+\frac{\epsilon-n}{2}}\right]
\]
for all dyadic cubes $Q$ and $R$, i.e.~that the operator
\[
  (b_Q)_Q \mapsto \bigg(\sum_{R} a_{QR}b_R\bigg)_Q
\]
is \emph{$\epsilon$-almost diagonal on $\dot{f}^{s,\tau}_{p,q}(\real^n)$} \cite[Definition 4.1]{YY2}. Now $\tau < \frac1p + \frac{\epsilon}{2n}$, so according to \cite[Theorem 4.1]{YY2}, the operator described above is bounded on $\dot{f}^{s,\tau}_{p,q}(\real^n)$. Also, by \cite[Theorem 3.1]{YY2}, the operator $u \mapsto \big(\langle u,\widetilde{\psi_Q}\rangle\big)_Q$ is bounded from $\dfp$ to $\dot{f}^{s,\tau}_{p,q}(\real^n)$. Combining these results with the estimates above yields
\[
  \|\widetilde{f}\|_{\mca^l_{N,m}\dfp} \lesssim \bigg\|\bigg(\sum_{R} a_{QR}|t_R|\bigg)_Q\bigg\|_{\dot{f}^{s,\tau}_{p,q}(\real^n)} \lesssim \left\|(t_Q)_Q\right\|_{\dot{f}^{s,\tau}_{p,q}(\real^n)} \lesssim \|f\|_{\dfp},
\]
i.e.~the mapping described above is a continuous embedding of $\dfp$ into $\mca^l_{N,m}\dfp$.

(ii) We have $\mca^l_{N,m}\dbp \subset \dbp$ in the same sense as above. Conversely, if a tempered distribution $f$ belongs to $\dbp$, again by the proof of \cite[Lemma 4.2]{YY2} and \cite[pp.~153--155]{FJ} it suffices to show that $\|\widetilde{f}\|_{\mca^l_{N,m}\dbp}$, where $\widetilde{f}$ is as in part (i), is controlled by a constant times $\|f\|_{\dbp}$. Using \eqref{eq:ad-estimate}, one can check that
\[
  \|\widetilde{f}\|_{\mca^l_{N,m}\dbp} \lesssim \bigg\|\bigg(\sum_{R} a_{QR}|t_R|\bigg)_Q\bigg\|_{\dot{b}^{s,\tau}_{p,q}(\real^n)},
\]
where $\dot{b}^{s,\tau}_{p,q}(\real^n)$ stands for the space of sequences $(b_Q)_Q$ of complex numbers, indexed by the dyadic cubes of $\real^n$, such that
\[
  \|(b_Q)_Q\|_{\dot{b}^{s,\tau}_{p,q}(\real^n)} = \sup_{P} \frac{1}{|P|^{\tau}}\Bigg(\sum_{j = j_P}^{\infty} \bigg[\sum_{Q \subset P,\;\ell(Q)=2^{-j}} \Big(|Q|^{-s/n-1/2 + 1/p} |b_Q| \Big)^q \bigg]^{p/q} \Bigg)^{1/p},
\]
where the supremum is taken over all dyadic cubes $P$ of $\real^n$ and $j_P = -\log_2\ell(P)$, is finite. Combining this estimate with \cite[Theorems 4.1 and 3.1]{YY2} as above then yields $\|\widetilde{f}\|_{\mca^l_{N,m}\dbp} \lesssim \|f\|_{\dbp}$.

(iii) This is contained in the arguments above.
\end{proof}

\section{Proof of Theorem \ref{th:pointwise}}\label{se:pointwise}

We now turn to the identification of the Haj\l asz-type spaces with the spaces defined through grand maximal functions. We shall need the following Sobolev-type embedding, which is a special case of \cite[Lemma 2.3]{KYZ2}.

\lem{le:embedding}
Let $s \in (0,1]$ and $0 < \epsilon < \epsilon' < s$. Then there exists a positive constant $C$ such that for all $x\in \real^n$, $k \in \integer$, measurable functions $u$ and $\overrightarrow{g} \in \D^s(u)$,
\begin{align*}
\inf_{c \in \complex} \dashint_{B(x,2^{-k})}|u(y) - c| dy \leq C 2^{-k\epsilon'}\sum_{j \geq k-2}2^{-j(s-\epsilon')}\left(\dashint_{B(x,2^{-k+1})}g_j(y)^{\frac{n}{n+\epsilon}}dy\right)^{\frac{n+\epsilon}{n}}.
\end{align*}
\elem

To talk about the identification of the spaces of measurable functions defined through Haj\l asz gradients and the spaces of tempered distributions defined through grand maximal functions, we need the following basic lemma. The techniques of the proof are similar to the ones employed in \cite[Theorem 1.1]{KYZ} and \cite[Theorem 3.2]{KYZ2}.

\lem{le:identifications}
(i) Let $u \in \dmp$ or $u \in \dnp$ with $s \in (0,1]$, $p \in (\frac{n}{n+s},\infty]$, $q \in (0,\infty]$ and $\tau \geq 0$. Then $u$ defines a tempered distribution in the sense that for all functions $\phi \in \sh(\real^n)$, $u\phi$ is integrable and
\[
  \Big| \int_{\real^n} u\phi \Big| \leq C(u) \|\phi\|_{\sh_{0,N}}
\]
for some integer $N$ depending on $n$, $s$, $p$ and $\tau$.

(ii) Suppose that $f\in\sdistr(\real^n)$ belongs to $\mca^l_{0,m}\dfp$ or $\mca^l_{0,m}\dbp$ with $s \in (0,\infty)$, $p \in (\frac{n}{n+s},\infty]$, $q \in (0,\infty]$, $\tau \in [0,\infty)$ and $l$, $m \in \nanu_0$. Then $f$ coincides with a locally integrable function in the sense that there exists a function $\tilde{f} \in L^1_{{\rm loc}}(\real^n)$ such that
\[
  \langle f, \phi \rangle = \int_{\real^n} \tilde{f}\phi
\]
for all $\phi \in \sh(\real^n)$ with compact support.
\elem

\begin{proof}
(i) Let $\|u\|$ denote either $\|u\|_{\dmp}$ or $\|u\|_{\dnp}$, whichever is finite. Fix $\epsilon$ and $\epsilon'$ such that $0 < \epsilon < \epsilon' < s$ and $\frac{n}{n+\epsilon} < p$. For any $x \in \real^n$, Lemma \ref{le:embedding} yields
\begin{align*}
  \inf_{c \in \complex} \dashint_{B(x,1)}|u(y) - c| dy & \lesssim \inf_{\overrightarrow{g} \in \D^s(u)}\sum_{j \geq -2}2^{-j(s-\epsilon')}\left(\dashint_{B(x,2)}g_j(y)^{\frac{n}{n+\epsilon}} dy\right)^{\frac{n+\epsilon}{n}}\\
& \lesssim \inf_{\overrightarrow{g} \in \D^s(u)} \sum_{j \geq -2}2^{-j(s-\epsilon')} \|g_j\|_{L^p(B(x,2))} \\
& \lesssim \sum_{j \geq -2} 2^{-j(s-\epsilon')} \|u\| < \infty,
\end{align*}
so $u$ is locally integrable.

Now write $N$ for the smallest integer strictly greater than $\max(n,n+s+n\tau - \frac{n}{p})$. For any $\phi \in \sh(\real^n)$ we then have
\begin{align*}
\int_{\real^n} |u(x)\phi(x)| dx & \lesssim |u|_{B(0,1)} \|\phi\|_{\sh_{0,N}} + \sum_{k \geq 1}\int_{B(0,2^k)\backslash B(0,2^{k-1})} \left|u(x) - u_{B(0,1)}\right||\phi(x)| dx \\
& \lesssim \|\phi\|_{\sh_{0,N}}\left( |u|_{B(0,1)} + \sum_{k \geq 1} 2^{-kN} \int_{B(0,2^k)} |u(x) - u_{B(0,1)}| dx\right) \\
& \approx \|\phi\|_{\sh_{0,N}}\left( |u|_{B(0,1)} + \sum_{k \geq 1} 2^{-k(N-n)} \dashint_{B(0,2^k)} |u(x) - u_{B(0,1)}| dx\right) \\
& \lesssim \|\phi\|_{\sh_{0,N}}\left( |u|_{B(0,1)} + \sum_{k \geq 1} 2^{-k(N-n)} \sum_{i=0}^{k}\dashint_{B(0,2^i)} |u(x) - u_{B(0,2^i)}| dx\right) \\
& \lesssim \|\phi\|_{\sh_{0,N}}\left( |u|_{B(0,1)} + \sum_{i \geq 0} 2^{-i(N-n)} \dashint_{B(0,2^i)} |u(x) - u_{B(0,2^i)}| dx\right).
\end{align*}
As above, the integral in the $i$th term of the latter sum can be estimated by
\begin{align*}
\dashint_{B(0,2^i)} |u(x) - u_{B(0,2^i)}| dx & \lesssim \inf_{\overrightarrow{g} \in \D^s(u)} 2^{i \epsilon'} \sum_{j \geq -i-2} 2^{-j(s-\epsilon')} \left(\dashint_{B(0,2^{i+1})}g_j(y)^{\frac{n}{n+\epsilon}} dy\right)^{\frac{n+\epsilon}{n}} \\
& \lesssim \inf_{\overrightarrow{g} \in \D^s(u)} 2^{i \epsilon'} \sum_{j \geq -i-2} 2^{-j(s-\epsilon')} |B(0,2^{i+1})|^{-\frac1p}\|g_j\|_{L^p(B(0,2^{i+1}))} \\
& \lesssim 2^{i(\epsilon'+n\tau-\frac{n}{p}}) \sum_{j \geq -i-2} 2^{-j(s-\epsilon')} \|u\| \\
& \lesssim 2^{i(s+n\tau - \frac{n}{p})} \|u\|.
\end{align*}
Thus,
\[
  \int_{\real^n} |u(x)\phi(x)| dx \lesssim \|\phi\|_{\sh_{0,N}}\left( |u|_{B(0,1)} + \sum_{i \geq 0} 2^{i(n+s+n\tau-\frac{n}{p}-N)} \|u\| \right),
\]
where the quantity inside the latter parentheses is finite because $n+s+n\tau - \frac{n}{p} - N < 0$.

(ii) Let $\|f\|$ denote either $\|f\|_{\mca^l_{0,m}\dfp}$ or $\|f\|_{\mca^l_{0,m}\dbp}$, whichever is finite. Fix a compactly supported $\vartheta \in \sh(\real^n)$ such that $\int_{\real^n} \vartheta = 1$. It is well known that
\[
  f = \vartheta * f + \sum_{j=0}^{\infty} \left(\vartheta_{j+1} - \vartheta_j\right)*f
\]
with convergence in $\sdistr(\real^n)$, so it suffices to show that
\[
  \sum_{j\geq 0} \left\| \left(\vartheta_{j+1}-\vartheta_j\right)*f\right\|_{L^1(B)} < \infty
\]
whenever $B \subset \real^n$ is a ball with radius $1$. For $p \geq 1$, this is almost immediate: $\vartheta_{1} - \vartheta$ is a constant multiple of an element of $\mca^l_{0,m}$, so
\[
  \sum_{j\geq 0} \left\| \left(\vartheta_{j+1}-\vartheta_j\right)*f\right\|_{L^1(B)} \lesssim \sum_{j\geq 0} 2^{-js}\big\| 2^{js}\left(\vartheta_{1}-\vartheta\right)_j*f\big\|_{L^p(B)} \lesssim \sum_{j\geq 0} 2^{-js} \|f\| < \infty.
\]
Suppose now that $\frac{n}{n+s} < p < 1$. If $x \in \real^n$ and $y \in B(x,2^{-j})$, $j \in \nanu_0$, then $\tilde{\phi}(z) := \phi(z - 2^j(x-y))$ is for all $\phi \in \mca^l_{0,m}$ a uniform constant multiple of some element of $\mca^l_{0,m}$, and we have $\tilde{\phi}_{j}(x) = \phi_j(y)$. Thus,
\begin{align*}
  \sup_{\phi\in\mca^l_{0,m}}|\phi_{j}*f(x)| & \lesssim \left(\dashint_{B(x,2^{-j})} \sup_{\phi\in\mca^l_{0,m}}|\phi_{j}*f(y)|^p dy\right)^{1/p} \\
& \lesssim 2^{\frac{jn}{p} - js} \Big\|2^{js} \sup_{\phi\in\mca^l_{0,m}}|\phi_{j}*f| \Big\|_{L^p(B(x,2^{-j}))} \lesssim 2^{\frac{jn}{p}-js-jn\tau}\|f\|.
\end{align*}
Using this estimate in an arbitrary ball $B$ of radius $1$ we thus get
\begin{align*}
\sum_{j \geq 0}\|(\vartheta_1 - \vartheta)_j*f\|_{L^1(B)}
& \lesssim \sum_{j \geq 0} \|(\vartheta_1-\vartheta)_j*f\|_{L^{\infty}(B)}^{1-p}\|(\vartheta_1 - \vartheta)_j*f\|_{L^p(B)}^p \\
& \lesssim \sum_{j \geq 0}\left(2^{\frac{jn}{p}-js-jn\tau}\|f\|\right)^{1-p}\left[2^{-js}\|f\|\right]^p \\
& \lesssim \sum_{j \geq 0} 2^{j \left(\frac{n}{p}(1-p) - s\right)} \|f\| < \infty,
\end{align*}
since $\frac{n}{p}(1-p) - s < (n+s)(1 - \frac{n}{n+s}) - s = 0$.
\end{proof}

We are now ready to give the proof of Theorem \ref{th:pointwise}. The methods are based on the case $p = \infty$ of the proof of \cite[Theorem 3.2]{KYZ2}.

\begin{proof}[Proof of Theorem \ref{th:pointwise}]
We shall first prove (i) and (ii) under the assumption that $s \in (0,1)$ and $m \geq n+1$.

(i) We start by establishing the embedding $\dmp \subset \mca^{\ell}_{0,m}\dfp$ with the assumption that $q < \infty$. To this direction, fix $\epsilon$ and $\epsilon'$ so that $0 < \epsilon < \epsilon' < s$ and $\frac{n}{n+\epsilon} < \min(p,q)$.

Let $u \in \dmp$, and choose $\overrightarrow{g} \in \D^s(u)$ so that
\[\sup_{x\in\real^n,\;\ell\in\integer} \frac{1}{|B(x,2^{-\ell})|^{\tau}}\left(\int_{B(x,2^{-\ell})} \bigg( \sum_{k \geq \ell}g_k(y)^q \bigg)^{p/q} dy\right)^{1/p} \leq 2\|u\|_{\dmp}. \]
Recall that by Lemma \ref{le:identifications} above, $u$ defines a tempered distribution. Since $m \geq n+1$, we have
\beqla{eq:h-g-estimate}
  \sup_{\phi \in \mca^l_{0,m}}|\phi_k * u(z)| \lesssim 2^{-k} \sum_{j \leq k} 2^{j(1-\epsilon')} \sum_{i \geq j-2} 2^{-i(s-\epsilon')}\left(\dashint_{B(z,2^{-j+1})} g_i(y)^{\frac{n}{n+\epsilon}} dy\right)^{\frac{n+\epsilon}{n}}
\eeq
for all $k \in \integer$ and $z \in \real^n$, where the implied constant does not depend on these two parameters; a similar estimate is established in \cite[pp.~15--16]{KYZ2}, but \eqref{eq:h-g-estimate} can also be deduced in a manner similar to the proof of part (i) of Lemma \ref{le:identifications}. We now consider a ball $B := B(x,2^{-\ell})$ with $\ell\in\integer$ and use the estimate \eqref{eq:h-g-estimate} for $z \in B$ and $k \geq \ell$. The terms of the sum with $j \leq \ell$ can be estimated as in the proof of Lemma \ref{le:identifications}:
\beqla{eq:h-g-subestimate1}
  \left(\dashint_{B(z,2^{-j+1})} g_i(y)^{\frac{n}{n+\epsilon}} dy\right)^{\frac{n+\epsilon}{n}} \lesssim 2^{\frac{jn}{p}}\|g_i\|_{L^p(B(z,2^{-j+2}))} \lesssim 2^{\frac{jn}{p} - jn\tau} \|u\|_{\dmp},
\eeq
and since $1 - s + \frac{n}{p} - n\tau > 0$, we get
\begin{align*}
  \sum_{j \leq \ell} 2^{j(1-\epsilon')} \sum_{i \geq j-2} 2^{-i(s-\epsilon')}\left(\dashint_{B(z,2^{-j+1})} g_i(y)^{\frac{n}{n+\epsilon}} dy\right)^{\frac{n+\epsilon}{n}}
& \lesssim \sum_{j \leq \ell} 2^{j - js + \frac{jn}{p} - jn\tau}\|u\|_{\dmp} \\
& \lesssim 2^{\ell(1 - s + \frac{n}{p} - n\tau)}\|u\|_{\dmp}.
\end{align*}
For the terms with $j > \ell$, we have $B(z,2^{-j+1}) \subset 2B$ for all $z \in B$, so that
\beqla{eq:h-g-subestimate2}
  \left(\dashint_{B(z,2^{-j+1})} g_i(y)^{\frac{n}{n+\epsilon}} dy\right)^{\frac{n+\epsilon}{n}} \leq \hlmax\left(g_i^{\frac{n}{n+\epsilon}} \chi_{2B}\right)(z)^{\frac{n+\epsilon}{n}},
\eeq
where $\hlmax$ is the Hardy-Littlewood maximal function. Combining these estimates, we have
\begin{align}
\; & \frac{1}{|B|^{\tau p}} \int_{B} \left(\sum_{k \geq \ell} \bigg[2^{ks}\sup_{\phi \in \mca^l_{0,m}}|\phi_k*u(z)|\bigg]^q\right)^\frac{p}{q} dz \notag \\
\lesssim\; & \frac{|B|}{|B|^{\tau p}} \left(\sum_{k\geq \ell} 2^{-k(1-s)q}\left(2^{\ell(1-s+\frac{n}{p}-n\tau)} \|u\|_{\dmp}\right)^q\right)^{\frac{p}{q}} \label{eq:h-g-subestimate3} \\
& \quad + \frac{1}{|B|^{\tau p}} \int_{B} \left(\sum_{k > \ell} \bigg[2^{-k(1-s)}\sum_{\ell < j \leq k} 2^{j(1-\epsilon')} \sum_{i \geq j-2} 2^{-i(s-\epsilon')} \hlmax\left(g_i^{\frac{n}{n+\epsilon}} \chi_{2B}\right)(z)^{\frac{n+\epsilon}{n}}\bigg]^q \right)^{\frac{p}{q}} dz.\notag
\end{align}
Since $1-s > 0$, the first quantity can be estimated easily:
\[
  \frac{|B|}{|B|^{\tau p}} \left(\sum_{k\geq \ell} 2^{-k(1-s)q}\left(2^{\ell(1-s+\frac{n}{p}-n\tau)} \|u\|_{\dmp}\right)^q\right)^{\frac{p}{q}} \lesssim \|u\|_{\dmp}^p.
\]
For the second quantity, we exchange the order of summation in the integrand as follows:
\begin{align*}
& \sum_{k > \ell} \left[2^{-k(1-s)}\sum_{\ell < j \leq k} 2^{j(1-\epsilon')} \sum_{i \geq j-2} 2^{-i(s-\epsilon')} \hlmax\left(g_i^{\frac{n}{n+\epsilon}} \chi_{2B}\right)(z)^{\frac{n+\epsilon}{n}}\right]^q \\
\approx & \sum_{k > \ell} \Bigg[2^{-k(1-s)} \sum_{i \geq \ell-1} 2^{-i(s-\epsilon')}\hlmax\left(g_i^{\frac{n}{n+\epsilon}} \chi_{2B}\right)(z)^{\frac{n+\epsilon}{n}} \sum_{\ell < j \leq \min(k,i+2)}2^{j(1-\epsilon')}\Bigg]^q \\
\approx & \sum_{k > \ell} \Bigg[2^{-k(1-s)} \sum_{\ell-1 \leq i \leq k-2} 2^{i(1-s)}\hlmax\left(g_i^{\frac{n}{n+\epsilon}} \chi_{2B}\right)(z)^{\frac{n+\epsilon}{n}} \\
& \quad \quad \quad + 2^{k(s-\epsilon')} \sum_{i \geq k-1} 2^{-i(s-\epsilon')}\hlmax\left(g_i^{\frac{n}{n+\epsilon}} \chi_{2B}\right)(z)^{\frac{n+\epsilon}{n}}  \Bigg]^q
\end{align*}
Now, using H\"older's inequality when $q > 1$ and the sub-additivity of $t \mapsto t^q$ when $\frac{n}{n+s} < q \leq 1$, the latter quantity can be estimated from above by a constant times
\begin{align*}
 & \sum_{k > \ell} \Bigg[2^{-k(1-s)(q\land1)} \sum_{\ell-1 \leq i \leq k-2} 2^{i(1-s)(q\land 1)}\hlmax\left(g_i^{\frac{n}{n+\epsilon}} \chi_{2B}\right)(z)^{\frac{n+\epsilon}{n}q} \\
& \quad \quad \quad + 2^{k(s-\epsilon')(q\land 1)} \sum_{i \geq k-1} 2^{-i(s-\epsilon')(q\land 1)}\hlmax\left(g_i^{\frac{n}{n+\epsilon}} \chi_{2B}\right)(z)^{\frac{n+\epsilon}{n}q}\Bigg]\\
\lesssim & \sum_{i \geq \ell -1} \hlmax\left(g_i^{\frac{n}{n+\epsilon}} \chi_{2B}\right)(z)^{\frac{n+\epsilon}{n}q}.
\end{align*}
Since $\frac{n+\epsilon}{n}p > 1$ and $\frac{n+\epsilon}{n}q > 1$, the Fefferman-Stein vector valued maximal inequality thus yields
\begin{align*}
\; & \frac{1}{|B|^{\tau p}} \int_{B} \left(\sum_{k > \ell} \bigg[2^{-k(1-s)}\sum_{\ell < j \leq k} 2^{j(1-\epsilon')} \sum_{i \geq j-2} 2^{-i(s-\epsilon')} \hlmax\left(g_i^{\frac{n}{n+\epsilon}} \chi_{2B}\right)(z)^{\frac{n+\epsilon}{n}}\bigg]^q \right)^{\frac{p}{q}} dz \\
\lesssim\; & \frac{1}{|B|^{\tau p}} \int_{B} \left(\sum_{i \geq \ell -1} \hlmax\left(g_i^{\frac{n}{n+\epsilon}} \chi_{2B}\right)(z)^{\frac{n+\epsilon}{n}q}\right)^\frac{\frac{n+\epsilon}{n}p}{\frac{n+\epsilon}{n}q} dz
\\
\lesssim\; &  \frac{1}{|B|^{\tau p}} \int_{2B} \bigg(\sum_{i\geq \ell - 1} g_i(z)^q\bigg)^{\frac{p}{q}} dz \lesssim \|u\|_{\dmp}^p.
\end{align*}
All in all, taking the supremum over all admissible $B$ in \eqref{eq:h-g-subestimate3} yields
\[
  \|u\|_{\mca^l_{0,m}\dfp} \lesssim \|u\|_{\dmp},
\]
i.e.~$\dmp \subset \mca^l_{0,m}\dfp$. The case $q = \infty$ can be handled in a similar but easier manner.

Now suppose that $f \in \mca^l_{0,m}\dfp$. According to Lemma \ref{le:identifications}, $f$ coincides with a locally integrable function, so fixing a compactly supported $\vartheta\in\sh(\real^n)$ such that $\int_{\real^n}\vartheta = 1$, we have $f = \lim_{j\to\infty} \vartheta_j*f$ pointwise almost everywhere. In particular, if $x$ and $y$ are distinct Lebesgue points of $f$ and $k\in\integer$ satisfies $2^{-k-1} \leq |x-y| < 2^{-k}$, we have
\[
  |f(x) - f(y)| \leq |\vartheta_k*f(x) - \vartheta_k*f(y)| + \sum_{j\geq k} \Big( |(\vartheta_1-\vartheta)_j*f(x)| + |(\vartheta_1-\vartheta)_j*f(y)|\Big),
\]
and since both $\vartheta_1 - \vartheta$ and $z \mapsto \vartheta(z) - \vartheta(z - 2^k[x-y])$ are constant multiples of elements of $\mca^l_{0,m}$, we get
\[
|f(x) - f(y)| \lesssim |x-y|^s (h_k(x) + h_k(y)),
\]
where
\[
  h_k(x) = 2^{ks} \sum_{j \geq k} \sup_{\phi \in \mca^l_{0,m}} |\phi_j*f(x)|
\]
for all $k \in \integer$ and $x\in \real$. Thus,
\[
  \|u\|_{\dmp} \lesssim \sup_{x\in\real^n,\;\ell \in \integer} \frac{1}{|B(x,2^{-\ell})|^{\tau}}\left(\int_{B(x,2^{-\ell})} \bigg( \sum_{k \geq \ell}h_k(z)^q \bigg)^{p/q} dz\right)^{1/p},
\]
and since for any $\ell \in \integer$ and $z \in \real^n$ we have the pointwise estimates
\begin{align*}
  \sum_{k \geq \ell}h_k(z)^q = & \sum_{k \geq \ell}2^{ksq} \bigg(\sum_{j \geq k} 2^{-js} \sup_{\phi \in \mca^l_{0,m}} 2^{js} |\phi_j*f(z)|\bigg)^q \\
\lesssim & \sum_{k \geq \ell}2^{ks(q\land 1)} \sum_{j \geq k} 2^{-js(q\land 1)}\sup_{\phi \in \mca^l_{0,m}} 2^{jsq}|\phi_j*f(x)|^q \\
\lesssim & \sum_{j \geq \ell} 2^{jsq} \sup_{\phi \in \mca^l_{0,m}} |\phi_j*f(z)|^q,
\end{align*}
we conclude that
\[
  \|u\|_{\dmp} \lesssim \|u\|_{\mca^l_{0,m}\dfp} ,
\]
i.e.~$\mca^l_{0,m}\dfp \subset \dmp$.

(ii) That $\dnp \subset \mca^l_{0,m}\dbp$ for $\frac{n}{n+s} < p < \infty$ can be established using the same estimates as above; here all values $q \in (0,\infty]$ are permitted because we only need the $L^{\frac{n+\epsilon}{n}p}$-boundedness of $\hlmax$, where $\epsilon \in (0,s)$ satisfies $p > \frac{n}{n+\epsilon}$, instead of the Fefferman-Stein maximal inequality. Let $u$ be a function in $\dnp$ and $\overrightarrow{g}\in\D^s(u)$ a fractional $s$-gradient that yields the quasinorm of $u$ up to a constant multiple of at most two. Letting $B := B(x,2^{-\ell})$ for arbitrary $x\in\real^n$ and $\ell\in\integer$, arguing as in \eqref{eq:h-g-estimate}, \eqref{eq:h-g-subestimate1} and \eqref{eq:h-g-subestimate2} we have
\begin{align*}
2^{ks}\sup_{\phi \in \mca^l_{0,m}}|\phi_k * u(z)| & \lesssim 2^{-k(1-s)}2^{\ell(1-s+\frac{n}{p}-n\tau)}\|u\|_{\dnp} \\
& \quad + 2^{-k(1-s)} \sum_{\ell < j \leq k} 2^{j(1-\epsilon')} \sum_{i\geq j-2} 2^{-i(s-\epsilon')}\hlmax\left(g_i^{\frac{n}{n+\epsilon}}\chi_{2B}\right)(z)^{\frac{n+\epsilon}{n}} \\
& \lesssim 2^{-k(1-s)}2^{\ell(1-s+\frac{n}{p}-n\tau)}\|u\|_{\dnp} \\
& \quad + 2^{-k(1-s)} \sum_{i \geq \ell -1}2^{-i(s-\epsilon')}2^{\min(k,i+2)(1-\epsilon')}\hlmax\left(g_i^{\frac{n}{n+\epsilon}}\chi_{2B}\right)(z)^{\frac{n+\epsilon}{n}},
\end{align*}
for all $k \geq \ell$ and $z \in B$, which further yields
\begin{align*}
\; & \big\|2^{ks}\sup_{\phi\in\mca^l_{0,m}}|\phi_k*f\big\|_{L^p(B)}^{p} \\
\lesssim\; & |B|^{\tau p}2^{(\ell-k)(1-s)p} \|u\|_{\dnp}^p + 2^{-k(1-s)p} \bigg\| \sum_{\ell - 1 \leq i \leq k-2} 2^{i(1-s)} \hlmax\left(g_i^{\frac{n}{n+\epsilon}}\chi_{2B}\right)^{\frac{n+\epsilon}{n}} \bigg\|_{L^p(B)}^p \\ 
& \;\quad\quad\quad\quad\quad\quad\quad\quad\quad\quad\quad+ 2^{k(s-\epsilon')p} \bigg\| \sum_{i \geq k - 1} 2^{-i(s-\epsilon')} \hlmax\left(g_i^{\frac{n}{n+\epsilon}}\chi_{2B}\right)^{\frac{n+\epsilon}{n}} \bigg\|_{L^p(B)}^p \\
\lesssim\; & |B|^{\tau p}2^{(\ell-k)(1-s)p} \|u\|_{\dnp}^p + 2^{-k(1-s)(p\land 1)} \sum_{\ell - 1 \leq i \leq k-2} 2^{i(1-s)(p\land 1)} \bigg\| \hlmax\left(g_i^{\frac{n}{n+\epsilon}}\chi_{2B}\right)^{\frac{n+\epsilon}{n}} \bigg\|_{L^p(B)}^p \\ 
& \;\quad\quad\quad\quad\quad\quad\quad\quad\quad\quad\quad+ 2^{k(s-\epsilon')(p\land 1)} \sum_{i \geq k -1} 2^{-i(s-\epsilon')(p\land 1)} \bigg\|\hlmax\left(g_i^{\frac{n}{n+\epsilon}}\chi_{2B}\right)^{\frac{n+\epsilon}{n}} \bigg\|_{L^p(B)}^p \\
\lesssim\; & |B|^{\tau p}2^{(\ell-k)(1-s)p} \|u\|_{\dnp}^p + 2^{-k(1-s)(p\land 1)} \sum_{\ell - 1 \leq i \leq k-2} 2^{i(1-s)(p\land 1)} \|g_i\|_{L^p(2B)}^p \\ 
& \;\quad\quad\quad\quad\quad\quad\quad\quad\quad\quad\quad+ 2^{k(s-\epsilon')(p\land 1)} \sum_{i \geq k -1} 2^{-i(s-\epsilon')(p\land 1)} \|g_i\|_{L^p(2B)}^p.
\end{align*}
Dividing out by $|B|^{\tau p} \approx |2B|^{\tau p}$, taking first the $\ell^{q/p}$ norm over $k \geq \ell$ and then taking the supremum over all admissible balls $B$ yields
\[
  \|u\|_{\mca^l_{0,m}\dbp} \lesssim \|u\|_{\dnp},
\]
i.e.~$\dnp \subset \mca^l_{0,m}\dbp$. The case $p = \infty$ can be handled in a similar but easier manner.

On the other hand, if $f \in \mca^l_{0,m}\dbp$ where $\frac{n}{n+s} < p < \infty$, then $f$ is a locally integrable function and $(h_k)_{k\in\integer}$ is a constant times an element of $\D^s(f)$, where the functions $h_k$ are as in part (i). For any ball $B := B(x,2^{-\ell})$, we have
\begin{align*}
|B|^{- \tau} \|h_k\|_{L^p(B)} \lesssim |B|^{-\tau} 2^{ks(p\land 1)(\frac1p\land 1)} \sum_{j \geq k} 2^{-js(p\land 1)(\frac1p\land 1)} \cdot 2^{js} \big\| \sup_{\phi\in\mca^l_{0,m}} |\phi_j * f| \big\|_{L^p(B)},
\end{align*}
so taking the $\ell^{q}$ norm over $k \geq \ell$ and then the supremum over all admissible balls $B$ yields
\[
  \|f\|_{\dnp} \lesssim \|f\|_{\mca^l_{0,m}\dbp},
\]
which means that $\mca^l_{0,m}\dbp \subset \dnp$. The case $p = \infty$ can be handled in a similar but easier manner.

The case with $s \in (0,1)$ and $m \geq n+1$ is thus proven. The case with $s = 1$ and $m \geq n+2$ can be proven using exactly the same argument, except for the fact that $1 - s$ is not positive. To remedy this, we use the assumption that $m \geq n+2$ to replace \eqref{eq:h-g-estimate} with the estimate
\[
  \sup_{\phi \in \mca^l_{0,m}}|\phi_k * u(z)| \lesssim 2^{-2k} \sum_{j \leq k} 2^{j(2-\epsilon')} \sum_{i \geq j-2} 2^{-i(s-\epsilon')}\left(\dashint_{B(z,2^{-j+1})} g_i(y)^{\frac{n}{n+\epsilon}} dy\right)^{\frac{n+\epsilon}{n}}.
\]
The above proof can thus be carried out by replacing $1-s$ with $2-s$ where necessary. We omit the details.
\end{proof}

\begin{acknowledgements}
The author would like to thank Eero Saksman and Yuan Zhou for reading the manuscript and making several valuable remarks. The author also wishes to thank IPAM at UCLA for supporting his participation in the program ``Interactions Between Analysis and Geometry'' in Spring 2013, during which part of this note was written. 
\end{acknowledgements}

\end{document}